\documentclass[11pt, reqno]{amsart}

\usepackage[usenames,dvipsnames]{pstricks}
\usepackage{epsfig}

\usepackage{color}
\date{\today}

\usepackage[latin1]{inputenc}
\usepackage[T1]{fontenc}
\usepackage{amsmath, amsthm}
\usepackage[english]{babel}
\usepackage{amssymb}
\usepackage{latexsym}
\usepackage{fancyhdr}
\usepackage{graphicx}
\usepackage[all]{xy}

\theoremstyle{definition}
\newtheorem{definicion}{{\bf Definition}}[section]

\newtheorem{teorema}[definicion]{{\bf Theorem}}

\newtheorem{corol}[definicion]{{\bf Corollary}}
\newtheorem{lema}[definicion]{{\bf Lemma}}

{{\noindent{\bf Proof: }}\newline }%
{$\hfill\Box$\vspace{0.25cm}} %

\newcommand{\ot}{\otimes}
\newcommand{\co}{\circ}

\begin{document}

\begin{center}
 {\huge{\bf  Partial and unified crossed products are weak }}

\vspace{0.25cm}

 {\huge{\bf crossed products}}
\end{center}

\ \\

{\bf  J.M. Fern\'andez Vilaboa$^{1}$, R. Gonz\'{a}lez
Rodr\'{\i}guez$^{2,\ast}$ and A.B. Rodr\'{\i}guez Raposo$^{3}$}

\ \\
\hspace{-0,5cm}$1$ Departamento de \'Alxebra, Universidad de
Santiago de Compostela,  E-15771 Santiago de Compostela, Spain
(e-mail: josemanuel.fernandez@usc.es)
\ \\
\hspace{-0,5cm}$2$ Departamento de Matem\'{a}tica Aplicada II,
Universidad de Vigo, Campus Universitario Lagoas-Marcosende,
E-36310 Vigo, Spain (e-mail: rgon@dma.uvigo.es)
\\
\hspace{-0,5cm}$3$ Departamento de Matem\'aticas, Universidade da Coru\~{n}a,  E-15071 A Coru\~{n}a, Spain
(e-mail: abraposo@edu.xunta.es)
\ \\
\hspace{-0,5cm}$\ast$ Corresponding author.

\begin{center}
{\bf Abstract}
\end{center}
{\small In \cite{mra-preunit} the notion of a weak crossed product
of an algebra by an object, both living in a  monoidal category was
presented. Unified crossed products defined in \cite{AM1} and
partial crossed products defined in \cite{partial} are crossed
product structures defined for a Hopf algebra and another object. In
this paper we prove that unified crossed products and partial
crossed products are particular instances of weak crossed products.
}

\vspace{0.5cm}



\section*{Introduction}
In \cite{nmra4}  an associative product on $A\otimes V$ was defined,
for an algebra $A$ and  an object $V$  both living in a strict monoidal category
$\mathcal C$ where every idempotent splits. This product was called the weak crossed product
of $A$ and $V$ and to obtain it we must consider crossed product systems, that is, two morphisms
$\psi_{V}^{A}:V\otimes A\rightarrow A\otimes V$ and $\sigma_{V}^{A}:V\otimes V\rightarrow A\otimes V$
satisfying some  twisted-like and cocycle-like conditions. Associated to these morphisms we define
an idempotent morphism $\nabla_{A\otimes V}:A\otimes V\rightarrow A\otimes V$.
The image of this idempotent, denoted by $A\times V$, inherits the associative product
from $A\otimes V$. In order to define a unit for $A\times V$, and hence to obtain an algebra
structure in this object,  in \cite{mra-preunit} we use the notion of  preunit introduced by
Caenepeel and De Groot in \cite{caengroot}. The theory presented in \cite{nmra4} and \cite{mra-preunit} contains the classical crossed products where $\nabla_{A\otimes V}=id_{A\ot V}$, for example the one defined by Brzezi\'nski in \cite{tb-crpr}, and also  other examples  with $\nabla_{A\otimes V}\neq id_{A\ot V}$ like the weak smash product given by Caenepeel and De Groot in \cite{caengroot}, the notion of weak wreath products that we can find in \cite{S} and the weak crossed products for weak bialgebras given in \cite{ana1} (see also \cite{mra-preunit}). Also,  B\"ohm showed in \cite{bohm} that a monad in the weak version of the Lack and Street's 2-category of monads in a 2-category is identical to a crossed product system in the sense of \cite{nmra4}.

Recently some new types of crossed products were presented in different settings
like for example partial crossed products and unified crossed products. The first one was
introduced by  Alves, Batista, Dokuchaev and  Paques \cite{partial} (see also the Batista's
presentation in the congress of Hopf Algebras and Tensor Categories that was held in Almer\'{\i}a
(Spain) July 4-8 (2011) http://www.ual.es/Congresos/hopf2010/charlas/batistalk.pdf) in order
to characterize cleft extensions of algebras in the partial setting. The notion of unified
crossed product was introduced by Agore and Militaru \cite{AM1} (see also \cite{AM2})  to describe and classify all Hopf algebras $E$ that factorize thorough $A$ and $H$ being $A$ a Hopf subalgebra of $E$, $H$ a subcoalgebra in $E$ with $1_{E}\in H$ and the multiplication $A\ot H\rightarrow E$ bijective.

In this paper we prove that partial and unified crossed products are weak crossed products.
The first one corresponds with a weak crossed product whose associated idempotent
is in general different of the identity  while in the second case the idempotent morphism is
the identity.

\section{The general theory of weak crossed products}

Throughout these notes  ${\mathcal C}$ denotes a strict
 monoidal category with tensor product $\otimes$ and base object $K$. Given objects
$A$, $B$, $D$ and a morphism
$f:B\rightarrow D$, we write $A\otimes f$ for $id_{A}\otimes f$
and $f\otimes A$ for $f\otimes id_{A}$. Also we assume  that
idempotents split, i.e., for every morphism
$\nabla_{Y}:Y\rightarrow Y$, such that $\nabla_{Y}=\nabla_{Y}\circ
\nabla_{Y}$, there exist an object $Z$ and morphisms
$i_{Y}:Z\rightarrow Y$ (injection) and $p_{Y}:Y\rightarrow Z$ (projection) satisfying
$\nabla_{Y}=i_{Y}\circ p_{Y}$ and $p_{Y}\circ i_{Y}=id_{Z}$.

An algebra (monoid) in ${\mathcal C}$ is a triple $A=(A, \eta_{A},
\mu_{A})$ where $A$ is an object in ${\mathcal C}$ and
 $\eta_{A}:K\rightarrow A$ (unit), $\mu_{A}:A\ot A
\rightarrow A$ (product) are morphisms in ${\mathcal C}$ such that
$\mu_{A}\co (A\ot \eta_{A})=id_{A}=\mu_{A}\co (\eta_{A}\ot A)$,
$\mu_{A}\co (A\ot \mu_{A})=\mu_{A}\co (\mu_{A}\ot A)$. Given two
algebras $A= (A, \eta_{A}, \mu_{A})$ and $B=(B, \eta_{B}, \mu_{B})$,
$f:A\rightarrow B$ is an algebra morphism if $\mu_{B}\co (f\ot
f)=f\co \mu_{A}$, $ f\co \eta_{A}= \eta_{B}$. Also, if ${\mathcal
C}$ is braided with braid $c$ and $A$, $B$ are algebras in
${\mathcal C}$, the object $A\otimes B$ is also an algebra in
 ${\mathcal C}$ where
$\eta_{A\otimes B}=\eta_{A}\otimes \eta_{B}$ and $\mu_{A\otimes
B}=(\mu_{A}\otimes \mu_{B})\circ (A\otimes c_{B,A}\otimes B).$

A coalgebra (comonoid) in ${\mathcal C}$ is a triple ${D} = (D,
\varepsilon_{D}, \delta_{D})$ where $D$ is an object in ${\mathcal
C}$ and $\varepsilon_{D}: D\rightarrow K$ (counit),
$\delta_{D}:D\rightarrow D\ot D$ (coproduct) are morphisms in
${\mathcal C}$ such that $(\varepsilon_{D}\ot D)\co \delta_{D}=
id_{D}=(D\ot \varepsilon_{D})\co \delta_{D}$, $(\delta_{D}\ot
D)\co \delta_{D}=
 (D\ot \delta_{D})\co \delta_{D}.$ If ${D} = (D, \varepsilon_{D},
 \delta_{D})$ and
${ E} = (E, \varepsilon_{E}, \delta_{E})$ are coalgebras,
$f:D\rightarrow E$ is a coalgebra morphism if $(f\ot f)\co
\delta_{D} =\delta_{E}\co f$, $\varepsilon_{E}\co f
=\varepsilon_{D}.$ When  ${\mathcal C}$ is braided with braid $c$
and $D$, $E$ are coalgebras in ${\mathcal C}$, $D\otimes E$ is a
coalgebra in ${\mathcal C}$ where $\varepsilon_{D\otimes
E}=\varepsilon_{D}\otimes \varepsilon_{E}$ and $\delta_{D\otimes
E}=(D\otimes c_{D,E}\otimes E)\circ( \delta_{D}\otimes \delta_{E}).$

If ${\mathcal C}$ is braided with braid $c$, we say that $H$ is a
bialgebra  in ${\mathcal C}$ if $(H, \eta_{H}, \mu_{H})$ is an
algebra, $(H, \varepsilon_{H}, \delta_{H})$ is a coalgebra and
$\varepsilon_{H}$ and $\delta_{H}$ are algebra morphisms
(equivalently $\eta_{H}$ and $\mu_{H}$ are coalgebra morphisms). If
moreover, there exists a morphism $\lambda_{H}:H\rightarrow H $
satisfying the identities
$$\mu_{H}\circ(H\otimes \lambda_{H})\circ \delta_{H}= \varepsilon_{H}\otimes
\eta_{H}= \mu_{H}\circ(\lambda_{H} \otimes H)\circ \delta_{H}$$ we
say that $H$ is a Hopf algebra.

 Let  $A$ be an algebra. The pair
$(M,\phi_{M})$ is a right $A$-module if $M$ is an object in
${\mathcal C}$ and $\phi_{M}:M\otimes A\rightarrow M$ is a morphism
in ${\mathcal C}$ satisfying $\phi_{M}\circ(M\otimes
\eta_{A})=id_{M}$, $\phi_{M}\circ (\phi_{M}\otimes A)=\phi_{M}\circ
(M\otimes \mu_{A})$. Given two right ${A}$-modules $(M,\phi_{M})$
and $(N,\phi_{N})$, $f:M\rightarrow N$ is a morphism of right
${A}$-modules if $\phi_{N}\circ (f\otimes A)=f\circ \phi_{M}$. In a
similar way we can define the notions of left $A$-module and
morphism of left $A$-modules. In this case we denote the left action
by $\varphi_{M}$.

In the first section of this note we develop the general theory of weak crossed products in a monoidal category $\mathcal C$ introduced in \cite{mra-preunit}.

Let $A$ be an algebra and $V$ be an object in
${\mathcal C}$. Suppose that there exists a morphism
$$\psi_{V}^{A}:V\otimes A\rightarrow A\otimes V$$  such that the following
equality holds
\begin{equation}\label{wmeas-wcp}
(\mu_A\otimes V)\circ (A\otimes \psi_{V}^{A})\circ
(\psi_{V}^{A}\otimes A) = \psi_{V}^{A}\circ (V\otimes \mu_A).
\end{equation}
 As a consequence of (\ref{wmeas-wcp}), the morphism $\nabla_{A\otimes V}:A\otimes V\rightarrow
A\otimes V$ defined by
\begin{equation}\label{idem-wcp}
\nabla_{A\otimes V} = (\mu_A\otimes V)\circ(A\otimes
\psi_{V}^{A})\circ (A\otimes V\otimes \eta_A)
\end{equation}
is  idempotent. Moreover, $\nabla_{A\otimes V}$
satisfies that $$\nabla_{A\otimes V}\circ (\mu_A\otimes V) =
(\mu_A\otimes V)\circ (A\otimes \nabla_{A\otimes V}),$$ that is,
$\nabla_{A\otimes V}$ is a left $A$-module morphism (see Lemma 3.1 of
\cite{mra-preunit}) for the regular action  $\varphi_{A\otimes V}=\mu_{A}\otimes V$. With $A\times V$, $i_{A\otimes V}:A\times V\rightarrow A\otimes V$ and $p_{A\otimes
V}:A\otimes V\rightarrow A\times V$ we denote the object, the injection and the projection
associated to the factorization of $\nabla_{A\otimes V}$.

From now on we consider quadruples $(A, V, \psi_{V}^{A}, \sigma_{V}^{A})$ where $A$
is an algebra, $V$ an object, $\psi_{V}^{A}:V\otimes A\rightarrow A\otimes V$ a morphism satisfiying (\ref{wmeas-wcp}) and
$\sigma_{V}^{A}:V\otimes V\rightarrow A\otimes V$  a morphism in ${\mathcal C}$.

We say that
$(A, V, \psi_{V}^{A},
\sigma_{V}^{A})$ satisfies the twisted condition if
\begin{equation}\label{twis-wcp}
(\mu_A\otimes V)\circ (A\otimes \psi_{V}^{A})\circ
(\sigma_{V}^{A}\otimes A) = (\mu_A\otimes V)\circ (A\otimes
\sigma_{V}^{A})\circ (\psi_{V}^{A}\otimes V)\circ (V\otimes
\psi_{V}^{A}).
\end{equation}
and   the  cocycle
condition holds if
\begin{equation}\label{cocy2-wcp}
(\mu_A\otimes V)\circ (A\otimes \sigma_{V}^{A}) \circ
(\sigma_{V}^{A}\otimes V) = (\mu_A\otimes V)\circ (A\otimes
\sigma_{V}^{A})\circ (\psi_{V}^{A}\otimes V)\circ
(V\otimes\sigma_{V}^{A}).
\end{equation}
By virtue of (\ref{twis-wcp}) and (\ref{cocy2-wcp}) we will consider from now
on, and without loss of generality, that
\begin{equation}
\label{idemp-sigma-inv}
\nabla_{A\otimes
V}\circ\sigma_{V}^{A} = \sigma_{V}^{A}
\end{equation}
holds for all quadruples $(A, V, \psi_{V}^{A}, \sigma_{V}^{A})$ {( see Proposition 3.7 of \cite{mra-preunit})}.

For $(A, V, \psi_{V}^{A}, \sigma_{V}^{A})$ define the product
\begin{equation}\label{prod-todo-wcp}
\mu_{A\otimes  V} = (\mu_A\otimes V)\circ (\mu_A\otimes
\sigma_{V}^{A})\circ (A\otimes \psi_{V}^{A}\otimes V)
\end{equation}
and let $\mu_{A\times V}$ be the product
\begin{equation}
\label{prod-wcp} \mu_{A\times V} = p_{A\otimes
V}\circ\mu_{A\otimes V}\circ (i_{A\otimes V}\otimes i_{A\otimes
V}).
\end{equation}

If the twisted and the cocycle conditions hold, the product $\mu_{A\otimes V}$ is associative and normalized with respect to $\nabla_{A\otimes
V}$ (i.e. $\nabla_{A\otimes
V}\circ \mu_{A\otimes V}=\mu_{A\otimes V}=\mu_{A\otimes V}\circ (\nabla_{A\otimes
V}\otimes \nabla_{A\otimes
V}$)). Due to this normality condition, $\mu_{A\times V}$ is associative as well (Propostion 3.7 of \cite{mra-preunit}). Hence we define:
\begin{definicion}\label{wcp-def}
If $(A, V, \psi_{V}^{A}, \sigma_{V}^{A})$  satisfies (\ref{twis-wcp}) and (\ref{cocy2-wcp})
we say that $(A\otimes V, \mu_{A\otimes V})$ is a weak
crossed product.
\end{definicion}
Our next aim is to endow $A\times V$ with a unit, and hence with an algebra structure. As $A\times V$ is given as an image of an idempotent, it seems reasonable to use a preunit on $A\otimes V$ to obtain a unit on $A\times V$. In general, if $A$ is an algebra, $V$ an object in ${\mathcal C}$ and $m_{A\otimes V}$ is an associative product defined in
$A\otimes V$ a preunit $\nu:K\rightarrow A\otimes V$ is a morphism satisfying
\begin{equation}
m_{A\otimes V}\circ (A\otimes V\otimes \nu)=m_{A\otimes V}\circ (\nu\otimes
A\otimes V)=m_{A\otimes V}\circ (A\otimes V\otimes (m_{A\otimes V}\circ
(\nu\otimes \nu))).
\end{equation}
Associated to a preunit we obtain an idempotent morphism
$$\nabla_{A\otimes
V}^{\nu}=m_{A\otimes V}\circ (A\otimes V\otimes \nu):A\otimes V\rightarrow
A\otimes V.$$
Take $A\times V$ the image of this idempotent, $p_{A\otimes V}^{\nu}$ the projection and $i_{A\otimes V}^{\nu}$ the injection. It is possible to endow $A\times V$ with an algebra structure whose product is $$m_{A\times V} = p_{A\otimes
V}^{\nu}\circ m_{A\otimes V}\circ (i_{A\otimes V}^{\nu}\otimes i_{A\otimes
V}^{\nu})$$
and whose unit is $\eta_{A\times V}=p_{A\otimes V}^{\nu}\circ \nu$ (see Proposition 2.5 of \cite{mra-preunit}). If moreover, $\mu_{A\otimes V}$ is left
$A$-linear for the actions $\varphi_{A\otimes V}=\mu_{A}\otimes V$, $\varphi_{A\otimes V\otimes
A\otimes V }=\varphi_{A\otimes V}\otimes  A\otimes V$ and normalized with respect to
$\nabla_{A\otimes
V}^{\nu}$,  the morphism
\begin{equation}
\label{beta-nu}
\beta_{\nu}:A\rightarrow A\otimes V,\; \beta_{\nu} =
(\mu_A\otimes V)\circ (A\otimes \nu)
\end{equation}
is
multiplicative and left $A$-linear for $\varphi_{A}=\mu_{A}$.

Although $\beta_{\nu}$ is not an algebra morphism, because
$A\otimes V$ is not an algebra, we have that $\beta_{\nu}\circ
\eta_A = \nu$, and thus the morphism $\bar{\beta_{\nu}} = p_{A\otimes
V}^{\nu}\circ\beta_{\nu}:A\rightarrow
A\times V$ is an algebra morphism.

In light of the considerations made in the last paragraphs, and
using the twisted and the cocycle conditions, in \cite{mra-preunit}
we characterize weak crossed products with a preunit, and moreover
we obtain an algebra structure on $A\times V$. These assertions are
a consequence of the following results proved in \cite{mra-preunit}.

\begin{teorema}
\label{thm1-wcp} Let $A$ be an algebra, $V$ an object and
$m_{A\otimes V}:A\otimes V\otimes A\otimes V\rightarrow A\otimes V$ a
morphism of left $A$-modules  for the actions $\varphi_{A\otimes V}=\mu_{A}\otimes
V$, $\varphi_{A\otimes V\otimes A\otimes V }=\varphi_{A\otimes V}\otimes  A\otimes V$.

Then the following statements are equivalent:
\begin{itemize}
\item[(i)] The product $m_{A\otimes V}$ is associative with preunit
$\nu$ and normalized with respect to $\nabla_{A\otimes V}^{\nu}.$

\item[(ii)] There exist morphisms $\psi_{V}^{A}:V\otimes A\rightarrow A\otimes
V$, $\sigma_{V}^{A}:V\otimes V\rightarrow A\otimes V$ and $\nu:k\rightarrow
A\otimes V$ such that if $\mu_{A\otimes V}$ is the product defined
in (\ref{prod-todo-wcp}), the pair $(A\otimes V, \mu_{A\otimes
V})$ is a weak crossed product with $m_{A\otimes V} =
\mu_{A\otimes V}$ satisfying:
    \begin{equation}\label{pre1-wcp}
    (\mu_A\otimes V)\circ (A\otimes \sigma_{V}^{A})\circ
    (\psi_{V}^{A}\otimes V)\circ (V\otimes \nu) =
    \nabla_{A\otimes V}\circ
    (\eta_A\otimes V),
    \end{equation}
    \begin{equation}\label{pre2-wcp}
    (\mu_A\otimes V)\circ (A\otimes \sigma_{V}^{A})\circ
    (\nu\otimes V) = \nabla_{A\otimes V}\circ (\eta_A\otimes V),
    \end{equation}
    \begin{equation}\label{pre3-wcp}
(\mu_A\otimes V)\circ (A\otimes \psi_{V}^{A})\circ (\nu\otimes A)
= \beta_{\nu},
\end{equation}
\end{itemize}
where $\beta_{\nu}$ is the morphism defined in
(\ref{beta-nu}). In this case $\nu$ is a preunit for $\mu_{A\otimes
V}$, the idempotent morphism of the weak crossed product
$\nabla_{A\otimes V}$ is the idempotent $\nabla_{A\otimes
V}^{\nu}$, and we say that the pair $(A\otimes V, \mu_{A\otimes V})$ is a
weak crossed product with preunit $\nu$.
\end{teorema}

\begin{corol}\label{corol-wcp}
If $(A\otimes V, \mu_{A\otimes V})$ is a weak crossed product with
preunit $\nu$, then $A\times V$ is an algebra with the product
defined in (\ref{prod-wcp}) and unit $\eta_{A\times V}=p_{A\otimes
V}\circ\nu$.
\end{corol}

\section{Partial crossed products are weak crossed products}

In this section we shall prove that the notion of partial crossed
product or crossed product by a partial action introduced in
\cite{partial}, for a Hopf algebra living in a category of
$k$-modules over an arbitrary unital commutative ring, is an example
of weak crossed product. To obtain this result we need that
${\mathcal C}$ be braided with braid $c$.

\begin{definicion}\label{twistedpartialaction}
Let $H$ be a Hopf algebra, $A$ an algebra and $\varphi_{A}:H\otimes
A\rightarrow A$, $\omega:H\otimes H\rightarrow A$ two morphisms in
${\mathcal C}$. The pair $(\varphi_{A}, \omega)$ is called a twisted
partial action of $H$ on $A$ if the following conditions hold:
\begin{itemize}
\item[(i)] $\varphi_{A}\circ (\eta_{H}\otimes A)=id_{A}$

\item[(ii)] $\varphi_{A}\circ (H\otimes \mu_{A})=\mu_{A}\circ (\varphi_{A}\otimes
\varphi_{A})\circ (H\otimes c_{H,A}\otimes A)\circ (\delta_{H}\otimes A\otimes A)$

\item[(iii)] $\mu_{A}\circ (\varphi_{A}\otimes \omega)\circ (H\otimes  c_{H,A}\otimes H)\circ (\delta_{H}\otimes ( (\varphi_{A}\otimes H)\circ (H\otimes c_{H,A})\circ (\delta_{H}\otimes A)))$
\item[ ]$=\mu_{A}\circ (A\otimes \varphi_{A})\circ (((\omega\otimes \mu_{H})\circ (H\otimes c_{H,H}\otimes H)\circ (\delta_{H}\otimes \delta_{H}))\otimes A)$

\item[(iv)] $\omega=\mu_{A}\circ (A\otimes \varphi_{A})\circ (((\omega\otimes \mu_{H})\circ (H\otimes c_{H,H}\otimes H)\circ (\delta_{H}\otimes \delta_{H}))\otimes \eta_{A}).$

\end{itemize}

If we write the previous identities in the monoidal category of $k$-modules over an arbitrary unital
commutative ring using elements and the Sweedler notation, we have the definition
of twisted partial action introduced in  \cite{partial}.

Moreover, if we define the morphisms
$$\psi_{H}^{A}:H\otimes A\rightarrow A\otimes H$$
by
\begin{equation}
\label{psi}
\psi_{H}^{A}=(\varphi_{A}\otimes H)\circ (H\otimes c_{H,A})\circ (\delta_{H}\otimes A)
\end{equation}
and
$$
\sigma_{H}^{A}:H\otimes H\rightarrow A\otimes H
$$
by
\begin{equation}
\label{sigma} \sigma_{H}^{A}=(\omega\otimes \mu_{H})\circ
\delta_{H\ot H}
\end{equation}
the equalities (ii), (iii) and (iv) of the previous definition can be rewritten in the  following form

\begin{itemize}

\item[(ii)] $\varphi_{A}\circ (H\otimes \mu_{A})=\mu_{A}\circ (A\otimes \varphi_{A})\circ (\psi_{H}^{A}\otimes A)$

\item[(iii)] $\mu_{A}\circ (A\otimes \omega)\circ (\psi_{H}^{A}\otimes H)\circ (H\otimes \psi_{H}^{A})
=\mu_{A}\circ (A\otimes \varphi_{A})\circ (\sigma_{H}^{A}\otimes A)$

\item[(iv)] $\omega=\mu_{A}\circ (A\otimes \varphi_{A})\circ (\sigma_{H}^{A}\otimes \eta_{A})$
\end{itemize}

\end{definicion}

The condition (iii) of the previous Definition can be seen as a
twisted condition in the partial setting.

Also, in \cite{partial} we
can find unit conditions
\begin{equation}
\label{unit}
\omega\circ (H\otimes \eta_{H})=\omega\circ (\eta_{H}\otimes H)=\varphi_{A}\circ (H\otimes \eta_{A})
\end{equation}
and a partial cocycle condition
\begin{equation}
\label{cocycle}
\mu_{A}\circ (\varphi_{A}\otimes \omega)\circ (H\otimes c_{H,A}\otimes H)\circ
(\delta_{H}\otimes ((\omega\otimes \mu_{H})\circ (H\otimes c_{H,H}\otimes H)\circ
(\delta_{H}\otimes \delta_{H})))
\end{equation}
$$=\mu_{A}\circ (A\otimes \omega)\circ (((\omega\otimes \mu_{H})\circ (H\otimes c_{H,H}\otimes H)
\circ (\delta_{H}\otimes \delta_{H}))\otimes H),$$
that admits the equivalent expression
\begin{equation}
\label{cocycle2}
\mu_{A}\circ (A\otimes \omega)\circ (\psi_{H}^{A}\otimes H)\circ (H\otimes \sigma_{H}^{A})=
\mu_{A}\circ (A\otimes \omega)\circ (\sigma_{H}^{A}\otimes H).
\end{equation}

\begin{lema}
Let $H$ be a Hopf algebra, $A$ an algebra and $\varphi_{A}:H\otimes
A\rightarrow A$, $\omega:H\otimes H\rightarrow A$ two morphisms. Let
$\psi_{H}^{A}$, $\sigma_{H}^{A}$ the morphisms defined in
(\ref{psi}) and (\ref{sigma}). Then the following equalities hold:
\begin{equation}
\label{sigdel} (\psi_{H}^{A}\ot H)\co (H\ot c_{H,A})\ot
(\delta_{H}\ot A)=(A\ot \delta_{H})\co \psi_{H}^{A} ,
\end{equation}
\begin{equation}
\label{sigdel1} (\sigma_{H}^{A}\ot \mu_{H})\co \delta_{H\ot H}=(A\ot
\delta_{H})\co \sigma_{H}^{A},
\end{equation}
\begin{equation}
\label{sigdel2} \varphi_{A}=(A\ot \varepsilon_{H})\co \psi_{H}^{A}.
\end{equation}
\begin{equation}
\label{sigdel3} \omega=(A\ot \varepsilon_{H})\co \sigma_{H}^{A},
\end{equation}
\end{lema}

\begin{proof} It is easy to prove that by the coassociativity of $\delta_{H}$ and the
naturality of $c$ we obtain (\ref{sigdel}). Using the same
properties we obtain (\ref{sigdel1}) because
\begin{itemize}
\item[ ] $\hspace{0.38cm}(\sigma_{H}^{A}\ot \mu_{H})\co \delta_{H\ot H} $
\item[ ] $=(\omega\ot \mu_{H}\ot \mu_{H})\co (H\ot c_{H,H}\ot c_{H,H}\ot H)\co
(H\ot c_{H,H}\ot H\ot H)\co (((\delta_{H}\ot H)\co \delta_{H})\ot $
\item[ ] $\hspace{0.38cm}
((\delta_{H}\ot H)\co \delta_{H}))$
\item[ ] $=(A\ot
\delta_{H})\co \sigma_{H}^{A}.$
\end{itemize}

Finally, the proof for (\ref{sigdel2}) is a consequence of the naturality of $c$ and the one for (\ref{sigdel3}) follows
from the fact that $\varepsilon_{H}$ is an algebra morphism.
\end{proof}

We have the following result.

\begin{teorema}
\label{idemp} Let $H$ be a Hopf algebra, $A$ an algebra and
$\varphi_{A}:H\otimes A\rightarrow A$, $\omega:H\otimes H\rightarrow
A$ two morphisms. If the condition (ii) of Definition
\ref{twistedpartialaction} is satisfied then the equality (\ref{wmeas-wcp})
holds for the morphism $\psi_{H}^{A}$  defined in
(\ref{psi}). As a consequence the morphism $\nabla_{A\otimes H}:A\otimes
V\rightarrow A\otimes H$ defined in (\ref{idem-wcp}) is  idempotent.
\end{teorema}

\begin{proof}
The equality (\ref{wmeas-wcp}) is fullfilled because
\begin{itemize}
\item[ ] $\hspace{0.38cm}(\mu_A\otimes H)\circ (A\otimes \psi_{H}^{A})\circ
(\psi_{H}^{A}\otimes A)    $
\item[ ] $=((\mu_{A}\co (\varphi_{A}\ot \varphi_{A}))\ot H)\co (H\ot c_{H,A}\ot c_{H,A})\co
(\delta_{H}\ot c_{H,A} \ot A)\co (\delta_{H}\ot A\ot A) $
\item[ ] $=((\varphi_{A}\circ (H\ot \mu_{A}))\ot H)\co (H\ot A\ot c_{H,A})\circ
(H\ot c_{H,A}\ot A)\co (\delta_{H}\ot A\ot A) $
\item[ ] $=\psi_{H}^{A}\co (H\ot \mu_{A}) $
\end{itemize}
where the first equality follows by the naturality of $c$ and the
coassociativity of $\delta_{H}$, the second one by (ii) of
Definition \ref{twistedpartialaction} and, finally, the last one by
the naturality of $c$.

Therefore by Lemma 3.1 of \cite{mra-preunit} we obtain that the
morphism  $\nabla_{A\otimes H}:A\otimes H\rightarrow A\otimes H$
defined in (\ref{idem-wcp}) is idempotent.

\end{proof}

Let ${\mathcal C}$ be  the monoidal category of $k$-modules over an
arbitrary unital commutative ring. If we write $\nabla_{A\otimes H}$
using the Sweedler notation and denoting $\varphi_{A}(h\otimes c)$
by $h.c$ we have
$$\nabla_{A\otimes H}(a\otimes h)=\sum a(h_{(1)}.{\bf 1}_{A})\otimes h_{(2)}.$$
Therefore, the  $k$-submodule generated by the elements $\sum
a(h_{(1)}.{\bf 1}_{A})\otimes h_{(2)}$ is $A\times H$ and then
$A\times H$ is the object  denoted in \cite{partial} by
$A\sharp_{(\varphi_{A},\omega)}H$.

Note that if  (\ref{unit}) holds we have:
\begin{equation}
\label{nabla-2} \nabla_{A\otimes H}=((\mu_{A}\co (A\ot \omega))\ot
H)\co (A\ot \eta_{H}\ot \delta_{H}).
\end{equation}

We have the following two Theorems involving the twisted an cocycle
conditions:

\begin{teorema}
\label{twistedTeo} Let $H$ be a Hopf algebra, $A$ an algebra and
$\varphi_{A}:H\otimes A\rightarrow A$, $\omega:H\otimes H\rightarrow
A$ two morphisms. Then the partial twisted condition given in (iii) of Definition
\ref{twistedpartialaction} holds if and only if the corresponding morphisms $\psi_{H}^{A}$ and
$\sigma_{H}^{A}$ defined in (\ref{psi}) and (\ref{sigma}) satisfy the twisted condition (\ref{twis-wcp}).
\end{teorema}

\begin{proof} If we assume that (\ref{twis-wcp}) holds for the morphisms
$\psi_{H}^{A}$ and $\sigma_{H}^{A}$ defined in (\ref{psi}) and
(\ref{sigma}), composing with $A\otimes \varepsilon_{H}$ in
(\ref{twis-wcp}) and using (\ref{sigdel2}) and (\ref{sigdel3}) we
obtain (iii) of Definition \ref{twistedpartialaction}. Conversely,
if  (iii) of Definition \ref{twistedpartialaction} holds we have the
equality (\ref{twis-wcp}) for the morphisms $\psi_{H}^{A}$ and
$\sigma_{H}^{A}$ defined in (\ref{psi}) and (\ref{sigma}). Indeed:

\begin{itemize}
\item[ ] $\hspace{0.38cm} (\mu_A\otimes H)\circ (A\otimes
\sigma_{H}^{A})\circ (\psi_{H}^{A}\otimes H)\circ (H\otimes
\psi_{H}^{A})   $
\item[ ] $=((\mu_{A}\circ (A\otimes \omega)\circ (\psi_{H}^{A}\otimes H)\circ
(H\otimes \psi_{H}^{A}))\otimes \mu_{H})\co (H\ot H\ot c_{H,A}\ot
H)\co (H\ot c_{H,H}\ot c_{H,A})\co $
\item[ ] $\hspace{0.38cm}(\delta_{H}\ot \delta_{H}\ot
A)$
\item[ ] $=((\mu_{A}\circ (A\otimes \varphi_{A})\circ (\sigma_{H}^{A}\otimes
A)) \otimes \mu_{H})\co (H\ot H\ot c_{H,A}\ot H)\co (H\ot c_{H,H}\ot
c_{H,A})\co (\delta_{H}\ot \delta_{H}\ot A)$
\item[ ] $=((\mu_{A}\co (A\ot \varphi_{A}))\ot H) \co (A\ot H\ot c_{H,A})\co
(((\sigma_{H}^{A}\ot \mu_{H})\co \delta_{H\ot H})\ot A)$
\item[ ] $=(\mu_A\otimes H)\circ (A\otimes \psi_{H}^{A})\circ
(\sigma_{H}^{A}\otimes A) $
\end{itemize}

In the last calculus the first equality follows by the naturality of
$c$ and the coassociativity of $\delta_{H}$, the second one by (iii)
of Definition \ref{twistedpartialaction}, the third  one by the
naturality of $c$ and the fourth one by (\ref{sigdel1}).
\end{proof}

\begin{teorema}
\label{cocycleTeo}
Let $H$ be a Hopf algebra, $A$ an algebra and
$\varphi_{A}:H\otimes A\rightarrow A$, $\omega:H\otimes H\rightarrow
A$ two morphisms. Then the partial cocycle condition given in (\ref{cocycle2}) holds if and only if the corresponding morphisms $\psi_{H}^{A}$ and
$\sigma_{H}^{A}$ defined in (\ref{psi}) and (\ref{sigma}) satisfy the cocycle condition (\ref{cocy2-wcp}).

\end{teorema}

\begin{proof} As in the previous Theorem, if we assume that (\ref{cocy2-wcp}) holds for the morphisms
$\psi_{H}^{A}$ and $\sigma_{H}^{A}$ defined in (\ref{psi}) and
(\ref{sigma}), composing with $A\otimes \varepsilon_{H}$  and using
(\ref{sigdel3}) we obtain (\ref{cocycle2}). Conversely, if
(\ref{cocycle2}) holds we have (\ref{cocy2-wcp}) for the morphisms
$\psi_{H}^{A}$ and $\sigma_{H}^{A}$ defined in (\ref{psi}) and
(\ref{sigma}). Indeed:

\begin{itemize}
\item[ ] $\hspace{0.38cm}(\mu_A\otimes H)\circ (A\otimes
\sigma_{H}^{A})\circ (\psi_{H}^{A}\otimes H)\circ
(H\otimes\sigma_{H}^{A})    $
\item[ ] $= ((\mu_{A}\circ (A\ot \omega))\ot\mu_{H})\co (A\ot H\ot c_{H,H}\ot H)\co (((A\ot \delta_{H})
\co \psi_{H}^{A})\ot H\ot H)\co (H\ot ((A\ot \delta_{H})\co \sigma_{H}^{A}))$
\item[ ] $=((\mu_{A}\circ (A\ot \omega))\ot\mu_{H})\co (A\ot H\ot c_{H,H}\ot H)\co (((\psi_{H}^{A}\ot H)
\co (H\ot c_{H,A})\co (\delta_{H}\ot A))\ot H\ot H)\co $
\item[ ] $\hspace{0.38cm} (H\ot ((\sigma_{H}^{A}\ot \mu_{H})\co
\delta_{H\ot H}))$
\item[ ] $=((\mu_{A}\circ (A\otimes \omega)\circ (\psi_{H}^{A}\otimes H)\circ
(H\otimes \sigma_{H}^{A}))\otimes \mu_{H})\co (H\ot H\ot c_{H,H}\ot \mu_{H})\co
(H\ot c_{H,H}\ot c_{H,H}\ot H)\co $
\item[ ] $\hspace{0.38cm}(\delta_{H} \ot \delta_{H} \ot \delta_{H})$
\item[ ] $=((\mu_{A}\circ (A\otimes \omega)\circ (\sigma_{H}^{A}\otimes H))\otimes \mu_{H})\co (H\ot H\ot c_{H,H}\ot \mu_{H})\co
(H\ot c_{H,H}\ot c_{H,H}\ot H)\co $
\item[ ] $\hspace{0.38cm}(\delta_{H} \ot \delta_{H} \ot \delta_{H})$
\item[ ] $=((\mu_{A}\circ (A\otimes \omega))\otimes \mu_{H})\co (A\ot H\ot c_{H,H}\ot H)\co
(((\sigma_{H}^{A}\ot \mu_{H})\co \delta_{H\ot H})\ot \delta_{H})$
\item[ ] $= (\mu_A\otimes H)\circ (A\otimes \sigma_{H}^{A}) \circ
(\sigma_{H}^{A}\otimes H)$
\end{itemize}

In the last calculus the first equality follows by definition of
$\sigma_{H}^{A}$, the second one by (\ref{sigdel}) and
(\ref{sigdel1}), the third one naturality of $c$ and the fourth one
by (\ref{cocycle2}). The fifth equality is a consequence of the
associativity of $\mu_{H}$ and the naturality of $c$ and, finally,
the last equality follows by (\ref{sigdel1}).
\end{proof}

As a consequence of this results, every twisted partial action
$(\varphi_{A},\omega)$ satisfying (\ref{cocycle2}) induces a weak
crossed product $(A\otimes H, \mu_{A\otimes H})$ where
$\psi_{H}^{A}:H\otimes A\rightarrow A\otimes H$ and
$\sigma_{H}^{A}:H\otimes H\rightarrow A\otimes H$ are the morphisms
defined in (\ref{psi}) and (\ref{sigma}).

The product defined in $A\otimes H$ in \cite{partial} by
$$(a\otimes h)(b\otimes l)=\sum a(h_{(1)}.b)\omega(h_{(2)}\otimes l_{(1)})\otimes h_{(3)}l_{(2)}$$
is $\mu_{A\otimes H}$ (i.e. the one induced by the quadruple $(A, H,
\psi_{H}^{A}, \sigma_{H}^{A})$) because in a monoidal setting the
previous equality can be written as
$$\mu_{A\otimes  H} = (\mu_A\otimes H)\circ (\mu_A\otimes
\sigma_{H}^{A})\circ (A\otimes \psi_{H}^{A}\otimes H).$$

Then, the product (called in \cite{partial} crossed product by a
twisted partial action or partial crossed product) induced in
$A\times H=A\sharp_{(\varphi_{A},\omega)}H$ is associative as well.

Moreover, we have the following: if $\nu=\eta_{A}\otimes \eta_{H}$,
by (ii) of Definition \ref{twistedpartialaction} and (\ref{unit}) we
obtain (\ref{pre1-wcp}). Similarly, by (\ref{unit}) we obtain
(\ref{pre2-wcp}) and finally (\ref{pre3-wcp}) follows by (i) of
Definition \ref{twistedpartialaction}. Therefore, if we assume unit
conditions we obtain that $\nu$ is a preunit (see Theorem
\ref{thm1-wcp}),  $(A\otimes H, \mu_{A\otimes H})$ is a weak crossed
product with preunit and then $A\times
H=A\sharp_{(\varphi_{A},\omega)}H$ is an algebra with unit
$$\eta_{A\times V}=p_{A\otimes V}^{\nu}\circ \nu=p_{A\otimes V}\circ \nu.$$

\section{Unified products are weak crossed products}

In the last  section of this paper we shall prove that the notion of
unified product introduced by Agore and Militaru in \cite{AM1} (see
also \cite{AM2}) is an example of a weak crossed product with
associate idempotent equal to the identity (note that in the
previous section $\nabla_{A\otimes H}\neq id_{A\otimes H}$). To
prove this assertion we also need that ${\mathcal C}$ be braided
with braid $c$.

The extension to the braided monoidal setting of the definition of
extending datum for a bialgebra $A$ introduced in \cite{AM1} is the
following:

\begin{definicion}
\label{ext-datum}
 Let $A$ be a bialgebra in ${\mathcal C}$. An
extending datum of $A$ is a system
$$\Omega(A)=(H, \phi_{H}:H\ot
A\rightarrow H, \varphi_{A}:H\ot A\rightarrow A, \tau:H\ot
H\rightarrow A)$$ where:
\begin{itemize}
\item[(i)] There exist morphisms $\eta_{H}:K\rightarrow H$,
$\mu_{H}:H\ot H\rightarrow H$, $\varepsilon_{H}:H\rightarrow K$ and
$\delta_{H}:H\rightarrow H\ot H$ such that
\begin{itemize}
\item[(i-1)] $(H,\varepsilon_{H},\delta_{H})$ is a coalgebra.
\item[(i-2)] $\delta_{H}\circ \eta_{H}=\eta_{H}\ot \eta_{H}.$
\item[(i-3)] $\mu_{H}\circ (\eta_{H}\ot H)=id_{H}=\mu_{H}\co (H\ot
\eta_{H})$.
\end{itemize}
\item[(ii)] $\phi_{H}$ is a coalgebra morphism i.e.
$$(\phi_{H}\ot \phi_{H})\co \delta_{H\ot A}=\delta_{H}\co
\phi_{H},\;\; \varepsilon_{H}\circ \phi_{H}=\varepsilon_{H}\ot
\varepsilon_{A}.$$

\item[(iii)] $\varphi_{A}$ is a coalgebra morphism i.e.
$$(\varphi_{A}\ot \varphi_{A})\co \delta_{H\ot A}=\delta_{A}\co
\varphi_{A},\;\; \varepsilon_{A}\circ \varphi_{A}=\varepsilon_{H}\ot
\varepsilon_{A}.$$

\item[(iv)] $\tau$ is a coalgebra morphism i.e.
$$(\tau\ot \tau)\co \delta_{H\ot H}=\delta_{A}\co
\tau,\;\; \varepsilon_{A}\circ \tau=\varepsilon_{H}\ot
\varepsilon_{H}.$$

\item[(v)] The following normalizing conditions hold
\begin{itemize}
\item[(v-1)] $\varphi_{A}\circ (H\otimes
\eta_{A})=\varepsilon_{H}\ot \eta_{A}$.
\item[(v-2)] $\varphi_{A}\circ (\eta_{H}\ot A)=id_{A}.$
\item[(v-3)] $\phi_{H}\circ (\eta_{H}\ot A)=\eta_{H}\ot
\varepsilon_{A}.$
\item[(v-4)] $\phi_{H}\circ (H\ot \eta_{A})=id_{H}.$
\item[(v-5)] $\tau\circ (H\ot \eta_{H})=\tau\co (\eta_{H}\co H)=\eta_{A}\ot \varepsilon_{H}.$
\end{itemize}
\end{itemize}

For an extending datum of $A$ we can define the morphisms

$$\psi_{H}^{A}:H\otimes A\rightarrow A\otimes H$$
by
\begin{equation}
\label{psi-ext} \psi_{H}^{A}=(\varphi_{A}\ot \phi_{H})\co
\delta_{H\ot A}
\end{equation}
and
$$
\sigma_{H}^{A}:H\otimes H\rightarrow A\otimes H
$$
by
\begin{equation}
\label{sigma-ext} \sigma_{H}^{A}=(\tau\otimes \mu_{H})\circ \delta_{H\ot
H}.
\end{equation}

Then, if we define the product
$$\mu_{A\otimes  H} = (\mu_A\otimes H)\circ (\mu_A\otimes
\sigma_{H}^{A})\circ (A\otimes \psi_{H}^{A}\otimes H),$$ in the
particular case of the category of $k$-modules for $k$ a field (or
an arbitrary unital commutative ring) we obtain that $\mu_{A\otimes
H}$ is the product $\bullet$ introduced in (16) of \cite{AM1} in the
following way:
$$(a\ot h)\bullet (c\ot g)=a(h_{(1)}\triangleright
c_{(1)})\tau((h_{(2)}\triangleleft c_{(2)})\ot g_{(1)})\ot
(h_{(3)}\triangleleft c_{(3)}).g_{(2)}$$  where we denoted the
product of two elements $x,y\in A$ by $xy$, $\phi_{H}(h\ot a)$ by
$h\triangleleft a$, $\varphi_{A}(h\ot a)$ by $h\triangleright a$ and
$\mu_{H}(h\ot g)=h.g$.

Following \cite{AM1}, if we denote by $A\ltimes H$ the $k$-module
together with the product $\bullet$, the object $A\ltimes H$ is
called the unified product of $A$ and $\Omega(A)$ if $A\ltimes H$ is
a bialgebra with the multiplication given by $\bullet$ and  by the tensor product of coalgebras. Agore and
Militaru proved in Theorem 2.4 of \cite{AM1} that $A\ltimes H$ is an
unified product if and only if $\delta_{H}$ and $\varepsilon_{H}$
are multiplicative maps, $(H,\phi_{H})$ is a rigth $A$-module
and the following equalities hold:

\begin{itemize}
\item[(BE1)] $(g\cdot h)\cdot l =
\bigl(g \triangleleft \tau(h_{(1)}\otimes l_{(1)})\bigl)\cdot (h_{(2)}\cdot l_{(2)})$

\item[(BE2)] $g \triangleright (ab) = (g_{(1)} \triangleright
a_{(1)})[(g_{(2)}\triangleleft a_{(2)})\triangleright b]$

\item[(BE3)] $(g\cdot h)
\triangleleft a = [g \triangleleft (h_{(1)} \triangleright a_{(1)})] \cdot (h_{(2)}
\triangleleft a_{(2)})$

\item[(BE4)] $[g_{(1)} \triangleright (h_{(1)}
\triangleright a_{(1)})]\tau\Bigl((g_{(2)} \triangleleft (h_{(2)} \triangleright
a_{(2)}))\otimes (h_{(3)} \triangleleft a_{(3)})\Bigl) = \tau(g_{(1)}\otimes h_{(1)})[(g_{(2)}
\cdot h_{(2)}) \triangleright a]$

\item[(BE5)] $\Bigl(g_{(1)} \triangleright
\tau(h_{(1)}\otimes  l_{(1)})\Bigl) \tau\Bigl((g_{(2)} \triangleleft \tau(h_{(2)}\otimes l_{(2)}))\otimes
(h_{(3)} \cdot l_{(3)})\Bigl) = \tau(g_{(1)}\otimes h_{(1)})\tau((g_{(2)} \cdot h_{(2)})\otimes l)$

\item[(BE6)] $g_{(1)} \triangleleft a_{(1)} \otimes g_{(2)} \triangleright a_{(2)} =
g_{(2)} \triangleleft a_{(2)} \otimes g_{(1)} \triangleright a_{(1)}$

\item[(BE7)]
$g_{(1)} \cdot h_{(1)} \otimes \tau(g_{(2)}\otimes h_{(2)}) = g_{(2)} \cdot h_{(2)} \otimes
\tau(g_{(1)}\otimes h_{(1)})$

\end{itemize}

for all $g$, $h$, $l \in H$ and $a$, $b \in
A$.

In the monoidal setting to say that $\delta_{H}$ and $\varepsilon_{H}$ are multiplicative morphisms is equivalent to
$$\delta_{H}\circ \mu_{H}=\mu_{H\ot H}\co (\delta_{H}\ot \delta_{H}),\;\;\; \varepsilon_{H}\circ \mu_{H}=\varepsilon_{H}\ot \varepsilon_{H}$$
and  $(H,\phi_{H})$ is a rigth $A$-module if
$$\phi_{H}\circ(H\otimes
\eta_{A})=id_{H},\;\;\;\phi_{H}\circ (\phi_{H}\otimes A)=\phi_{H}\circ
(M\otimes \mu_{A}).$$
On the other hand, the equalities (BE1) to (BE7)  can be rewritten  using the morphisms $\psi_{H}^{A}$ and $\sigma_{H}^{A}$ defined in (\ref{psi-ext}) and (\ref{sigma-ext}) in the  following form:

\begin{itemize}

\item[(BE1)] $\mu_{H}\co (\mu_{H}\ot H)=\mu_{H}\co (\phi_{H}\ot H)\co (H\ot \sigma_{H}^{A})$

\item[(BE2)] $\varphi_{A}\co (H\ot \mu_{A})=\mu_{A}\co (A\ot \varphi_{A})\co (\psi_{H}^{A}\ot A)$

\item[(BE3)] $\phi_{H}\co (\mu_{H}\ot A)=\mu_{H}\co (\phi_{H}\ot H)\co (H\ot \psi_{H}^{A})$

\item[(BE4)] $\mu_{A}\circ (A\otimes \tau)\circ (\psi_{H}^{A}\otimes H)\circ (H\otimes \psi_{H}^{A})
=\mu_{A}\circ (A\otimes \varphi_{A})\circ (\sigma_{H}^{A}\otimes A)$

\item[(BE5)] $\mu_{A}\circ (A\otimes \tau)\circ (\psi_{H}^{A}\otimes H)\circ (H\otimes \sigma_{H}^{A})=
\mu_{A}\circ (A\otimes \tau)\circ (\sigma_{H}^{A}\otimes H)$

\item[(BE6)] $c_{A,H}\circ \psi_{H}^{A}=(\phi_{H}\ot \varphi_{A})\co \delta_{H\ot A}$

\item[(BE7)] $c_{A,H}\circ \sigma_{H}^{A}=(\mu_{H}\ot \tau)\co \delta_{H\ot H}$

\end{itemize}

\end{definicion}

Note that (BE4) and (BE5) are similar to (iii) of Definition \ref{twistedpartialaction} and (\ref{cocycle2}) respectively. Then we call them the unified twisted condition and the unified cocycle condition respectively.

\begin{lema}
Let $\Omega(A)$ be an extending datum of a bialgebra $A$. Let
$\psi_{H}^{A}$, $\sigma_{H}^{A}$ be the morphisms defined in
(\ref{psi-ext}) and (\ref{sigma-ext}). Then the following equalities hold:
\begin{equation}
\label{sigdel-ext} (\psi_{H}^{A}\ot \phi_{H})\co \delta_{H\ot A}=(A\ot \delta_{H})\co \psi_{H}^{A} ,
\end{equation}
\begin{equation}
\label{sigdel-extnew} (\varphi_{A}\ot \psi_{H}^{A})\co \delta_{H\ot A}=(\delta_{A}\ot H)\co \psi_{H}^{A} ,
\end{equation}
\begin{equation}
\label{sigdel-extnew1}
(\delta_{A}\ot H)\co \sigma_{H}^{A}=(\tau\ot \sigma_{H}^{A})\co \delta_{H\ot H},
\end{equation}
\begin{equation}
\label{sigdel1-ext} \varphi_{A}=(A\ot \varepsilon_{H})\co \psi_{H}^{A},
\end{equation}
\begin{equation}
\label{sigdel1-extphi} \phi_{H}=(\varepsilon_{A}\ot H)\co \psi_{H}^{A},
\end{equation}
\begin{equation}
\label{sigdel1-extsigma} \mu_{H}=(\varepsilon_{A}\ot H)\co \sigma_{H}^{A}.
\end{equation}
Moreover, if $\delta_{H}$ is a multiplicative morphism, the equality
\begin{equation}
\label{sigdel2-ext} (\sigma_{H}^{A}\ot \mu_{H})\co \delta_{H\ot H}=(A\ot
\delta_{H})\co \sigma_{H}^{A}
\end{equation}
holds, and if $\varepsilon_{H}$ is a multiplicative morphism we have the identity
\begin{equation}
\label{sigdel3-ext} \tau=(A\ot \varepsilon_{H})\co \sigma_{H}^{A},
\end{equation}
\end{lema}

\begin{proof} It is easy to prove that by the condition of coalgebra morphism for $\phi_{H}$, the coassociativity of $\delta_{H}$ and the
naturality of $c$  we obtain (\ref{sigdel-ext}). In a similar way, using the condition of of coalgebra morphism for $\varphi_{A}$, the coassociativity of $\delta_{H}$ and the
naturality of $c$ we prove (\ref{sigdel-extnew}). The proof for (\ref{sigdel-extnew1}) is equal to the previous ones using  the condition  of coalgebra morphism for $\tau$. By $\varepsilon_{A}\co \varphi_{A}=\varepsilon_{H}\ot \varepsilon_{A}$ we obtain (\ref{sigdel1-extphi}) and (\ref{sigdel1-extsigma}) follows by $\varepsilon_{A}\co \tau=\varepsilon_{H}\ot \varepsilon_{H}$.
The equality (\ref{sigdel1-ext}) follows directly from the condition of coalgebra morphism for $\phi_{H}$. Finally,
the proofs for (\ref{sigdel2-ext}) and (\ref{sigdel3-ext}) are similar with the ones used for (\ref{sigdel1}) and (\ref{sigdel3}).
\end{proof}

\begin{teorema}
\label{newidemp} Let $\Omega(A)$ be an extending datum of a
bialgebra $A$ satisfying (BE2) and such that $(H,\phi_{H})$ is a
rigth $A$-module. Then the equality (\ref{wmeas-wcp}) holds for the
morphism $\psi_{H}^{A}$  defined in (\ref{psi-ext}).
\end{teorema}

\begin{proof}
Let $\psi_{H}^{A}$ be the morphism defined in (\ref{psi-ext}). Then

\begin{itemize}
\item[ ] $\hspace{0.38cm}(\mu_A\otimes H)\circ (A\otimes
\psi_{H}^{A})\circ (\psi_{H}^{A}\otimes H)$
\item[ ] $=\mu_{A}\circ (A\ot \varphi_{A}\ot \phi_{H})\co (A\ot H\ot c_{H,A}\ot A)\co
(((A\ot \delta_{H})\co \psi_{H}^{A})\ot \delta_{A})  $
\item[ ] $=((\mu_{A}\co (A\ot \varphi_{A})\co (\psi_{H}^{A}\ot A))\ot \phi_{H})\co
(H\ot A\ot (c_{H,A}\co (\phi_{H}\ot A))\ot A)\co (\delta_{H\ot A}\ot \delta_{A}) $
\item[ ] $=((\varphi_{A}\co (H\ot \mu_{A}))\ot (\phi_{H}\co (\phi_{H}\ot A)))\co
(H\ot A\ot c_{H,A}\ot A\ot A)\co (H\ot c_{H,A}\ot c_{A,A}\ot A)\co $
\item[ ] $\hspace{0.38cm} (\delta_{H} \ot \delta_{A}\ot \delta_{A}) $
\item[ ] $= (\varphi_{A}\ot \phi_{H})\co (H\ot c_{H,A}\ot A)\co (\delta_{H}\ot
((\mu_{A}\ot \mu_{A})\co \delta_{A\ot A})) $
\item[ ] $= \psi_{H}^{A}\co (H\ot \mu_{A})$
\end{itemize}
where the first equality follows by definition, the second one by (\ref{sigdel-ext}),
the third one by (BE2), the fourth one by the condition of right $A$-module for $H$ and
the naturality of $c$ and finally, the in the last one we used that $A$ is a bialgebra.
\end{proof}

Then, as a consequence of the previous Theorem, we have that for any extending datum of a bialgebra
$A$ satisfying (BE2) and such that  $(H,\phi_{H})$ is a right $A$-module,
$$(A, H, \psi_{H}^{A}=(\varphi_{A}\ot \phi_{H})\co
\delta_{H\ot A}, \sigma_{H}^{A}=(\tau\ot \mu_{H})\co
\delta_{H\ot H})$$
 is a quadruple as the ones used to define the notion of weak crossed product. In this case the
 associate idempotent $\nabla_{A\ot H}$ is the identity morphism because, by the normalizing
 conditions ((v) of Definition \ref{ext-datum}), we obtain
$$\nabla_{A\ot H}=((\mu_{A}\co (A\ot (\varphi_{A}\co (H\ot \eta_{A}))))\ot (\phi_{H}\co (H\ot
\eta_{A})))\co (A\ot \delta_{H})=id_{A\ot H}$$
Therefore, in this case $A\times H=A\ot H$.

\begin{lema}
\label{apoyo1} Let $\Omega(A)$ be an extending datum of a bialgebra
$A$ satisfying (BE6). Then the equality
\begin{equation}
(A\ot (c_{H,H}\co (\phi_{H}\ot H)))\co (c_{H,A}\ot A\ot H)\co (H\ot ((\delta_{A}\ot H)\co
\psi_{H}^{A}))=
\end{equation}
$$(\psi_{H}^{A}\ot \phi_{H})\co (H\ot c_{H,A}\ot \varphi_{A})\co (c_{H,H}\ot c_{H,A}\ot A)\co
(H\ot \delta_{H}\ot \delta_{A})$$
 holds for the morphism $\psi_{H}^{A}$  defined in (\ref{psi-ext}).
\end{lema}

\begin{proof} We have
\begin{itemize}
\item[ ] $\hspace{0.38cm}(A\ot (c_{H,H}\co (\phi_{H}\ot H)))\co (c_{H,A}\ot A\ot H)\co
(H\ot ((\delta_{A}\ot H)\co
\psi_{H}^{A})) $
\item[ ] $= (H\ot (c_{H,H}\co (\phi_{H}\ot H)))\co (c_{H,A}\ot A\ot H)\co (H\ot
((\varphi_{A}\ot \psi_{H}^{A})\co \delta_{H\ot A})) $
\item[ ] $= (A\ot H\ot \phi_{H})\co (A\ot c_{H,H}\ot A)\co
(((\varphi_{A}\ot H)\co (H\ot c_{H,A}))\ot
(c_{A,H}\co\psi_{H}^{A}))\co (c_{H,H}\ot c_{H,A}\ot A)\co  $
\item[ ] $\hspace{0.38cm} (H\ot \delta_{H}\ot \delta_{A})$
\item[ ] $= (A\ot H\ot \phi_{H})\co (A\ot c_{H,H}\ot A)\co
(((\varphi_{A}\ot H)\co (H\ot c_{H,A}))\ot ((\phi_{H}\ot
\varphi_{A})\co \delta_{H\ot A}))\co $
\item[ ] $\hspace{0.38cm}  (c_{H,H}\ot c_{H,A}\ot A)\co  (H\ot \delta_{H}\ot \delta_{A})$
\item[ ] $=(\psi_{H}^{A}\ot \phi_{H})\co (H\ot c_{H,A}\ot \varphi_{A})\co (c_{H,H}\ot c_{H,A}\ot A)\co
(H\ot \delta_{H}\ot \delta_{A}). $
\end{itemize}
where the first equality follows by (\ref{sigdel-extnew1}), the
second one by the naturality of the braiding, the third one by (BE6)
and the last one by the coassociativity of $\delta_{H}$,
$\delta_{A}$  and the naturality of $c$.

\end{proof}

\begin{lema}
\label{apoyo2} Let $\Omega(A)$ be an extending datum of a bialgebra
$A$ satisfying (BE7). Then the equality
\begin{equation} (A\ot (c_{H,H}\co (\phi_{H}\ot H)))\co (c_{H,A}\ot
A\ot H)\co (H\ot ((\delta_{A}\ot H)\co \sigma_{H}^{A}))=
\end{equation}
$$(\sigma_{H}^{A}\ot \phi_{H})\co (H\ot c_{H,A}\ot \tau)\co (c_{H,H}\ot c_{H,H}\ot H)\co
(H\ot \delta_{H}\ot \delta_{H})$$
 holds for the morphism $\sigma_{H}^{A}$  defined in (\ref{sigma-ext}).
\end{lema}

\begin{proof} We have
\begin{itemize}
\item[ ] $\hspace{0.38cm}(A\ot (c_{H,H}\co (\phi_{H}\ot H)))\co (c_{H,A}\ot
A\ot H)\co (H\ot ((\delta_{A}\ot H)\co \sigma_{H}^{A})) $

\item[ ] $= (H\ot (c_{H,H}\co (\phi\ot H)))\co (c_{H,A}\ot A\ot H)\co (H\ot
((\tau\ot \sigma_{H}^{A})\co \delta_{H\ot H})) $

\item[ ] $= (A\ot H\ot \phi_{H})\co (A\ot c_{H,H}\ot A)\co
(((\tau\ot H)\co (H\ot c_{H,H}))\ot (c_{A,H}\co\sigma_{H}^{A}))\co
(c_{H,H}\ot c_{H,H}\ot H)\co  $

\item[ ] $\hspace{0.38cm} (H\ot \delta_{H}\ot \delta_{H})$

\item[ ] $= (A\ot H\ot \phi_{H})\co (A\ot c_{H,H}\ot A)\co
(((\tau\ot H)\co (H\ot c_{H,H}))\ot ((\mu_{H}\ot \tau)\co
\delta_{H\ot H}))\co $

\item[ ] $\hspace{0.38cm}  (c_{H,H}\ot c_{H,H}\ot H)\co  (H\ot \delta_{H}\ot \delta_{H})$

\item[ ] $=(\sigma_{H}^{A}\ot \phi_{H})\co (H\ot c_{H,H}\ot \tau)\co (c_{H,H}\ot c_{H,H}\ot H)\co
(H\ot \delta_{H}\ot \delta_{H}). $
\end{itemize}
where the first equality follows by (\ref{sigdel-extnew1}), the
second one by the naturality of the braiding, the third one by (BE7)
and the last one by the coassociativity of $\delta_{H}$,
  and the naturality of $c$.

\end{proof}

\begin{teorema}
\label{newtwisted1} Let $\Omega(A)$ be an extending datum of a
bialgebra $A$ such that $\varepsilon_{H}$ is multiplicative. Then if
the twisted condition (\ref{twis-wcp}) holds for the morphisms
$\psi_{H}^{A}$, $\sigma_{H}^{A}$ defined in (\ref{psi-ext}) and
(\ref{sigma-ext}), the unified twisted condition (BE4) holds.
\end{teorema}

\begin{proof} The proof is a direct consequence of
(\ref{sigdel1-ext}) and (\ref{sigdel3-ext}).
\end{proof}

\begin{teorema}
\label{newtwisted2} Let $\Omega(A)$ be an extending datum of a
bialgebra $A$ satisfying (BE3), the unified twisted condition (BE4), (BE6) and such that
$\delta_{H}$ is multiplicative. Then the twisted condition
(\ref{twis-wcp}) holds for the morphisms $\psi_{H}^{A}$,
$\sigma_{H}^{A}$ defined in (\ref{psi-ext}) and (\ref{sigma-ext}).
\end{teorema}

\begin{proof} We have
\begin{itemize}
\item[ ] $\hspace{0.38cm} (\mu_A\otimes H)\circ (A\otimes
\sigma_{H}^{A})\circ (\psi_{H}^{A}\otimes H)\circ (H\otimes
\psi_{H}^{A})$
\item[ ] $= (\mu_{A}\ot H)\co (A\ot \tau\ot \mu_{H})\co (A\ot H\ot c_{H,H}\ot H)\co
(((A\ot \delta_{H})\co \psi_{H}^{A})\ot H\ot H)\co $
\item[ ] $\hspace{0.38cm}(H\ot ((A\ot \delta_{H})\co \psi_{H}^{A}))) $
\item[ ] $=((\mu_{A}\circ (A\ot \tau))\ot \mu_{H})\circ (\psi_{H}^{A}\ot H\ot H\ot H)\co   $
\item[ ] $\hspace{0.38cm}(H\ot ((A\ot (c_{H,H}\co (\phi_{H}\ot H)))\co (c_{H,A}\ot A\ot H)\co
(H\ot ((\delta_{A}\ot H)\co
\psi_{H}^{A}))) \ot \phi_{H})\co $
\item[ ] $\hspace{0.38cm} (\delta_{H}\ot \delta_{H\ot A})$
\item[ ] $=((\mu_{A}\circ (A\ot \tau))\ot \mu_{H})\circ (\psi_{H}^{A}\ot H\ot H\ot H)\co   $
\item[ ] $\hspace{0.38cm}(H\ot ((\psi_{H}^{A}\ot \phi_{H})\co (H\ot c_{H,A}\ot \varphi_{A})\co (c_{H,H}\ot c_{H,A}\ot A)\co
(H\ot \delta_{H}\ot \delta_{A})) \ot \phi_{H})\co $
\item[ ] $\hspace{0.38cm} (\delta_{H}\ot \delta_{H\ot A})$
\item[ ] $= (\mu_{A}\ot \mu_{H})\co (A\ot \varphi_{A}\ot \phi_{H}\ot H)\co
(\sigma_{H}^{A}\ot c_{H,A}\ot \varphi_{A}\ot H)\co (H\ot c_{H,H}\ot
c_{H,A}\ot A\ot H)\co  $
\item[ ] $\hspace{0.38cm} (\delta_{H}\ot \delta_{H}\ot \delta_{A}\ot
\phi_{H})\co (H\ot \delta_{H\ot A})$
\item[ ] $= ((\mu_{A}\co (A\ot \varphi_{A})\co ( \sigma_{H}^{A}\ot A))\ot
(\mu_{H}\co (\phi_{H}\ot H)\co
(H\ot \psi_{H}^{A})))\co (H\ot H\ot c_{H,A}\ot H\ot A)\co $
\item[ ] $\hspace{0.38cm}(H\ot c_{H,H}\ot c_{H,A}\ot A)\co
(\delta_{H}\ot \delta_{H}\ot \delta_{A})$
\item[ ] $=((\mu_{A}\co (A\ot \varphi_{A})\co ( \sigma_{H}^{A}\ot A))\ot
(\phi_{H}\co (\mu_{H}\ot A)))\co (H\ot H\ot c_{H,A}\ot H\ot A)\co $
\item[ ] $\hspace{0.38cm}
(H\ot c_{H,H}\ot c_{H,A}\ot A)\co (\delta_{H}\ot \delta_{H}\ot
\delta_{A})$
\item[ ] $= ((\mu_{A}\co (A\ot \varphi_{A}))\ot \phi_{H})\co (A\ot H\ot c_{H,A}\ot A)\co
(((\sigma_{H}^{A}\ot \mu_{H})\co \delta_{H\ot H})\ot \delta_{A}) $
\item[ ] $= (\mu_{A}\ot H)\co (A\ot ((\varphi_{A}\ot \phi_{H})\co \delta_{H\ot A}))\co
(\sigma_{H}^{A}\ot A) $
\item[ ] $= (\mu_A\otimes H)\circ (A\otimes \psi_{H}^{A})\circ
(\sigma_{H}^{A}\otimes A) $
\end{itemize}
where the first equality follows by definition, the second one by
(\ref{sigdel-ext}), the third one by Lemma \ref{apoyo1} and the
fourth one by (BE4). The fifth equality is a consequence of the
coassociativity of  $\delta_{H}$, $\delta_{A}$  and the naturality
of $c$. In the sixth one we used (BE3) and the seventh follows by
the naturality of $c$. Finally, the eighth follows by
(\ref{sigdel2-ext}) ($\delta_{H}$ is multiplicative) and the last
one by definition.

\end{proof}

\begin{teorema}
\label{newcocycle1} Let $\Omega(A)$ be an extending datum of a
bialgebra $A$ such that $\varepsilon_{H}$ is multiplicative. Then if
the cocycle condition (\ref{cocy2-wcp}) holds for the morphisms
$\psi_{H}^{A}$, $\sigma_{H}^{A}$ defined in (\ref{psi-ext}) and
(\ref{sigma-ext}), the unified cocycle condition (BE5) holds.
\end{teorema}

\begin{proof} The proof is a direct consequence of
 (\ref{sigdel3-ext}).
\end{proof}

\begin{teorema}
\label{newcocycle2} Let $\Omega(A)$ be an extending datum of a
bialgebra $A$ satisfying (BE1), the unified cocycle condition  (BE5), (BE7) and such that
$\delta_{H}$ is multiplicative. Then the cocycle condition
(\ref{cocy2-wcp}) holds for the morphisms $\psi_{H}^{A}$,
$\sigma_{H}^{A}$ defined in (\ref{psi-ext}) and (\ref{sigma-ext}).
\end{teorema}

\begin{proof}
 We have
\begin{itemize}
\item[ ] $\hspace{0.38cm} (\mu_A\otimes H)\circ (A\otimes
\sigma_{H}^{A})\circ (\psi_{H}^{A}\otimes H)\circ (H\otimes
\sigma_{H}^{A})$

\item[ ] $= (\mu_{A}\ot H)\co (A\ot \tau\ot \mu_{H})\co (A\ot H\ot c_{H,H}\ot H)\co
(((A\ot \delta_{H})\co \psi_{H}^{A})\ot H\ot H)\co $
\item[ ] $\hspace{0.38cm}(H\ot ((A\ot \delta_{H})\co \sigma_{H}^{A}))) $

\item[ ] $=((\mu_{A}\circ (A\ot \tau))\ot \mu_{H})\circ (\psi_{H}^{A}\ot H\ot H\ot H)\co   $
\item[ ] $\hspace{0.38cm}(H\ot ((A\ot (c_{H,H}\co (\phi_{H}\ot H)))\co (c_{H,A}\ot A\ot H)\co
(H\ot ((\delta_{A}\ot H)\co \sigma_{H}^{A}))) \ot \mu_{H})\co $
\item[ ] $\hspace{0.38cm} (\delta_{H}\ot \delta_{H\ot H})$

\item[ ] $=((\mu_{A}\circ (A\ot \tau))\ot \mu_{H})\circ (\psi_{H}^{A}\ot H\ot H\ot H)\co   $
\item[ ] $\hspace{0.38cm}(H\ot ((\sigma_{H}^{A}\ot \phi_{H})\co (H\ot c_{H,A}\ot \tau)
\co (c_{H,H}\ot c_{H,H}\ot H)\co (H\ot \delta_{H}\ot \delta_{H}))
\ot \mu_{H})\co $
\item[ ] $\hspace{0.38cm} (\delta_{H}\ot \delta_{H\ot H})$

\item[ ] $= (\mu_{A}\ot \mu_{H})\co (A\ot \tau\ot \phi_{H}\ot H)\co
(\sigma_{H}^{A}\ot c_{H,H}\ot \tau\ot H)\co (H\ot c_{H,H}\ot
c_{H,H}\ot H\ot H)\co  $
\item[ ] $\hspace{0.38cm} (\delta_{H}\ot \delta_{H}\ot \delta_{H}\ot
\mu_{H})\co (H\ot \delta_{H\ot H})$

\item[ ] $= ((\mu_{A}\co (A\ot \tau)\co ( \sigma_{H}^{A}\ot H))\ot (\mu_{H}\co
(\phi_{H}\ot H)\co (H\ot \sigma_{H}^{A})))\co (H\ot H\ot c_{H,H}\ot
H\ot H)\co $
\item[ ] $\hspace{0.38cm}(H\ot c_{H,H}\ot c_{H,H}\ot H)\co
(\delta_{H}\ot \delta_{H}\ot \delta_{H})$

\item[ ] $=((\mu_{A}\co (A\ot \tau)\co ( \sigma_{H}^{A}\ot A))\ot
(\mu_{H}\co (\mu_{H}\ot H)))\co (H\ot H\ot c_{H,H}\ot H\ot H)\co $
\item[ ] $\hspace{0.38cm}
(H\ot c_{H,H}\ot c_{H,H}\ot H)\co (\delta_{H}\ot \delta_{H}\ot
\delta_{H})$

\item[ ] $= ((\mu_{A}\co (A\ot \tau))\ot \mu_{H})\co (A\ot H\ot c_{H,H}\ot H)\co
(((\sigma_{H}^{A}\ot \mu_{H})\co \delta_{H\ot H})\ot \delta_{H}) $

\item[ ] $= (\mu_{A}\ot H)\co (A\ot ((\tau\ot \mu_{H})\co \delta_{H\ot H}))\co
(\sigma_{H}^{A}\ot A) $

\item[ ] $= (\mu_A\otimes H)\circ (A\otimes \sigma_{H}^{A})\circ
(\sigma_{H}^{A}\otimes H) $

\end{itemize}
where the first equality follows by definition, the second one by
(\ref{sigdel-ext}) and (\ref{sigdel2-ext}), the third one by Lemma
\ref{apoyo2} and the fourth one by (BE5). The fifth equality is a
consequence of the coassociativity of  $\delta_{H}$ and the
naturality of $c$. In the sixth one we used (BE1) and the seventh
follows by the naturality of $c$. Finally, the eighth follows by
(\ref{sigdel2-ext}) ($\delta_{H}$ is multiplicative) and the last
one by definition.
\end{proof}

Then, as a consequence of Theorems \ref{newidemp}, \ref{newtwisted1},
\ref{newtwisted2}, \ref{newcocycle1} and \ref{newcocycle2},  if
$\Omega(A)$ is an extending datum of a bialgebra $A$ such that
$\delta_{H}$ and $\varepsilon_{H}$ are multiplicative morphisms,
$(H,\phi_{H})$ is a rigth $A$-module and the equalities (BE1) to
(BE7) hold, for the morphisms $\psi_{H}^{A}$, $\sigma_{H}^{A}$
defined in (\ref{psi-ext}) and (\ref{sigma-ext}), we obtain that
$(A,H,\psi_{H}^{A},\sigma_{H}^{A})$ is a quadruple as the ones
considered in \cite{mra-preunit} and the pair $(A\ot H,\mu_{A\otimes
H})$ with
$$\mu_{A\otimes  H} = (\mu_A\otimes H)\circ (\mu_A\otimes
\sigma_{H}^{A})\circ (A\otimes \psi_{H}^{A}\otimes H)$$ is a weak
crossed product that in the case of ${\mathcal C=k-}Mod $ is the
unified product defined in \cite{AM1}. In this case the preunit is a
unit because $\nabla_{A\ot H}=id_{A\ot H}$ and, by the normalizing
conditions of the Definition of extending datum, the equalities
$$\mu_{A\ot H}\co (\eta_{A}\ot \eta_{H}\ot A\ot H)=id_{A\ot
H}=\mu_{A\ot H}\co(A\ot H\ot \eta_{A}\ot \eta_{H})$$ hold.
Therefore, $(A\ot H, \eta_{A\ot H}=\eta_{A}\ot \eta_{H}, \mu_{A\ot
H})$ is an algebra in ${\mathcal C}$.
\section*{Acknowledgements}
The authors were supported by  Ministerio de Ciencia e Innovaci\'on
(Projects: MTM2010-15634, MTM2009-14464-C02-01), Xunta de Galicia
(Project: INCITE09207215PR) and by FEDER.


\begin{thebibliography}{99}

\bibitem{AM1}  A. Agore, G. Militaru:  Extending structures II:
The quantum version. arXiv:1011.2174v3 (2011)

\bibitem{AM2}  A. Agore, G. Militaru:
Unified products and split extensions of Hopf algebras.
arXiv:1105.1474v1 (2011)

\bibitem{nmra4} J.N. Alonso Álvarez, J.M. Fernández Vilaboa,
R. González Rodríguez, A.B. Rodríguez Raposo:  Crosssed products
in weak contexts. Appl. Cat. Structures 18, 231-258 (2010)

\bibitem{bohm} G. B\"ohm: The weak theory of monads,
 Adv. Math. 225, 1-32 (2010)

\bibitem{tb-crpr} T. Brzezi\'nski: Crossed products by a
coalgebra. Comm. in Algebra 25, 3551-3575  (1997)

\bibitem{caengroot} S. Caenepeel, E. De Groot: Modules over
weak entwining structures.  Contemp. Math. 267, 31-54 (2000)

\bibitem{mra-preunit}  J.M. Fern\'andez Vilaboa,
R. Gonz\'alez Rodr\'{\i}guez and A.B. Rodr\'{\i}guez Raposo:
Preunits and weak crossed products. J. of Pure Appl. Algebra 213,
2244-2261 (2009)

\bibitem{mra-proj} J.M. Fern\'andez Vilaboa, R. Gonz\'alez Rodr\'{\i}guez, A.B. Rodr\'{\i}guez Raposo:
Weak crossed biproducts and weak projections. arXiv:0906.1693v2
(2009)

\bibitem{partial}  M. Muniz S. Alves, E. Batista, M. Dokuchaev, A. Paques:
Twisted partial actions of Hopf algebras. preprint (2011)

\bibitem{ana1} A.B. Rodríguez Raposo: Crossed products for
weak Hopf algebras. Comm. in Algebra 37, 2274-2289  (2009)

\bibitem{S} R. Street: Weak distributive laws. Theory and Applications of Categories 22, 313-320 (2009)


\end{thebibliography}
\end{document}